\documentclass[12pt,fleqn]{article}
\usepackage{amsmath,amsthm,amssymb,amscd}
\usepackage{graphicx}
\usepackage{graphics}
\usepackage{verbatim}

\usepackage{mathrsfs}
\usepackage{color}
\makeatletter

\newcommand{\Rmnum}[1]{\expandafter\@slowromancap\romannumeral #1@}

\newtheorem{theorem}{Theorem}[section]

\newtheorem{lemma}[theorem]{Lemma}

\newtheorem{corollary}[theorem]{Corollary}

\makeatother

\begin{document}
\title{Remarks on planar edge-chromatic critical graphs}

\author{Ligang Jin\thanks{supported by Deutsche Forschungsgemeinschaft (DFG) grant STE 792/2-1;
Paderborn Institute for Advanced Studies in
		Computer Science and Engineering,
		Paderborn University,
		Warburger Str. 100,
		33102 Paderborn,
		Germany; ligang@mail.upb.de}, 
Yingli Kang\thanks{Fellow of the International Graduate School "Dynamic Intelligent Systems";
Paderborn Institute for Advanced Studies in
		Computer Science and Engineering,
		Paderborn University,
		Warburger Str. 100,
		33102 Paderborn,
		Germany; yingli@mail.upb.de}, 
Eckhard Steffen\thanks{
		Paderborn Institute for Advanced Studies in
		Computer Science and Engineering,
		Paderborn University,
		Warburger Str. 100,
		33102 Paderborn,
		Germany; es@upb.de}}
\date{}

\maketitle

\abstract
{The only open case of Vizing's conjecture that every planar graph with $\Delta\geq 6$ is a class 1 graph is $\Delta = 6$.
We give a short proof of the following statement:
there is no 6-critical plane graph $G$, such that every vertex of $G$ is incident to at most three 3-faces.
A stronger statement without restriction to critical graphs is stated in \cite{Wang_Xu_2013}.
However, the proof given there works only for critical graphs.
Furthermore, we show that every 5-critical plane graph has a 3-face which is adjacent to a $k$-face $(k\in \{3,4\})$.

For $\Delta = 5$ our result gives insights into the structure of planar $5$-critical graphs, and the result for $\Delta=6$
gives support for the truth of Vizing's planar graph conjecture.
}


\section{Introduction}

We consider finite simple graphs $G$ with vertex set $V(G)$ and edge set $E(G)$.
The vertex-degree of $v \in V(G)$ is denoted by $d_G(v)$, and $\Delta(G)$ denotes the maximum vertex-degree of $G$.
If it is clear from the context, then $\Delta$ is frequently used.
A graph is planar if it is embeddable into the Euclidean plane. A plane graph $(G,\Sigma)$ is a planar graph $G$ together with an embedding $\Sigma$ of $G$ into the Euclidean plane.
If $(G,\Sigma)$ is a plane graph, then $F(G)$ denotes the set of faces of $(G,\Sigma)$.
The degree $d_{(G,\Sigma)} (f) $ of a face $f$ is the length of its facial circuit. A face $f$ is a $k$-face if $d_{G}(f)=k$,
and it is a $k^+$-face if $d_{G}(f) \geq k$.

The edge-chromatic number $\chi'(G)$ of a graph $G$ is the minimum $k$ such that $G$ admits a proper $k$-edge-coloring.
Vizing \cite{Vizing_1964} proved that $\chi'(G) \in \{\Delta(G), \Delta(G)+1\}$. If $\chi'(G)= \Delta(G)$, then $G$ is a class 1 graph, and it is a
class 2 graph otherwise. A class 2 graph $H$ is $k$-critical, if $\Delta(H)=k$ and $\chi'(H') < \chi'(H)$ for every proper subgraph $H'$ of $H$.

Vizing \cite{Vizing_1964} showed for each $k \in \{2,3,4,5\}$ that there is a planar class 2 graph $G$ with $\Delta(G) = k$. He proved that
every planar graph with $\Delta\geq 8$ is a class 1 graph, and conjectured that every planar graph with $\Delta \in \{6,7\}$ is a class 1 graph. Vizing's conjecture is
proved for planar graphs with $\Delta=7$ by Gr\"unewald \cite{Gruenewald_2000},  Sanders, Zhao \cite{Sanders_Zhao_2001}, and Zhang \cite{Zhang_2000} independently.
It is still open for the case $\Delta=6$. The paper provides short proofs for the following statements.

\begin{theorem} \label{th_main}
There is no 6-critical plane graph $(G, \Sigma)$, such that every vertex of $G$ is incident to at most three 3-faces.
\end{theorem}

If Vizing's conjecture is not true, then every 6-critical graph has the following property.

\begin{corollary}
Let $(G,\Sigma)$ be a plane graph. If $G$ is $6$-critical, then there is a vertex of $G$ which is incident to at least four $3$-faces.
\end{corollary}

\begin{theorem} \label{th_main 1}
Let $(G,\Sigma)$ be a plane graph. If $G$ is $5$-critical, then $(G,\Sigma)$ has a $3$-face which is adjacent to a $3$-face or to a $4$-face.
\end{theorem}

\begin{figure}[h]
  \centering
  \includegraphics[width=4.5cm]{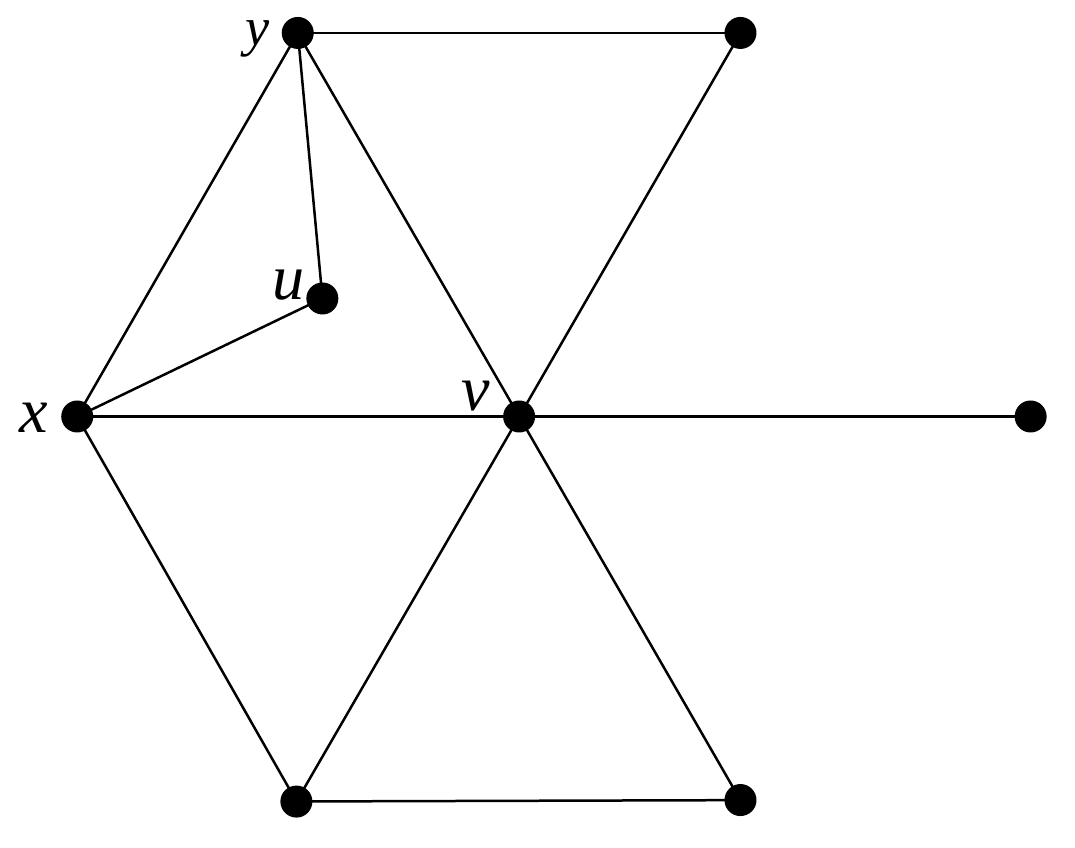}\\
  \caption{an example}\label{mistake_example}
\end{figure}

A significant longer proof of Theorem \ref{th_main} is given in \cite{Wang_Xu_2013}, but the statement is formulated for plane graphs.
However, the proof works for critical graphs only. The assumption that a minimal counterexample is critical is wrong. It might be that a subgraph of this minimal counterexample $G$ does not fulfill the pre-condition of the statement. For example, if $G$ has a triangle $[vxyv]$ and a bivalent vertex $u$ such that $u$ is the unique vertex inside $[vxyv]$ and $u$ is adjacent to $x$ and $y$, then the removal of $u$ increases the number of 3-faces containing $v$ (see Figure \ref{mistake_example}).

\section{Proofs of Theorems \ref{th_main} and \ref{th_main 1}}

We will use the following two lemmas.

\begin{lemma} [\cite{Luo_Miao_Zhao_2009}] \label{lem_m_n}
If $G$ is a $6$-critical graph, then $|E(G)|\geq \frac{1}{2}(5|V(G)|+3)$.
\end{lemma}

\begin{lemma} [\cite{Woodall_2008}] \label{lem_m_n 1}
If $G$ is a $5$-critical graph, then $|E(G)| \geq \frac{15}{7}|V(G)|$.
\end{lemma}

\subsection*{Proof of Theorem \ref{th_main}}

Suppose to the contrary that there is a counterexample to the statement. Then there is a $6$-critical graph $G$ which has an embedding $\Sigma$
such that every $v \in V(G)$ is incident to at most three 3-faces.
With Euler's formula and Lemma \ref{lem_m_n} we deduce $\sum_{f\in F(G)}(d_G(f)-4) = 2|E(G)|-4|F(G)| = 2|E(G)|-4(|E(G)|+2-|V(G)|) \leq -|V(G)|-11$.
Therefore,
$ |V(G)| + \sum_{f\in F(G)}(d_{G}(f)-4) \leq -11$.

Give initial charge 1 to each $v \in V(G)$ and $d_{G}(f)-4$ to each $f \in F(G)$. Discharge the elements of $V(G) \cup F(G)$ according to the following rule:\\
\textbf{R1}: Every vertex sends $\frac{1}{3}$ to its incident 3-faces.

The rule only moves the charge around and does not affect the sum. Furthermore, the finial charge of every vertex and face is at least 0.
Therefore, $0\leq\sum_{v\in V(G)}1 + \sum_{f\in F(G)}(d_{G}(f)-4) = |V(G)| + \sum_{f\in F(G)}(d_{G}(f)-4) \leq -11$, a contradiction.

\subsection*{Proof of Theorem \ref{th_main 1}}

Suppose to the contrary that there is a counterexample to the statement. Then there is a $5$-critical graph $G$ which has an embedding $\Sigma$
such that every 3-face is adjacent to $5^+$-faces only. Hence, every vertex of $G$ is incident to at most two 3-faces, and every vertex which is incident to a 3-face is also incident
to a $5^+$-face.
By Lemma \ref{lem_m_n 1}, we have $\sum_{f\in F(G)}(d_{G}(f)-4) \leq -\frac{2}{7}|V(G)|-8$. Therefore,
$\frac{2}{7}|V(G)| + \sum_{f\in F(G)}(d_{G}(f)-4) \leq -8$.

Give initial charge of $\frac{2}{7}$ to each vertex and $d_G(f)-4$ to each face of $G$. Discharge the elements of $V(G) \cup F(G)$ according to the following rules:\\
\textbf{R1}: Every vertex sends $\frac{1}{3}$ to its incident 3-faces.\newline
\textbf{R2}: Every $5^+$-face sends $\frac{d_G(f)-4}{d_G(f)}$ to its incident vertices.

Denote the finial charge by $ch^*$. Rules R1 and R2 imply that  $ch^*(f)\geq 0$ for every $f \in F(G)$.
Let $n \leq 2$ and $v$ be a vertex which is incident to $n$ 3-faces.
If $n=0$, then $ch^*(v)\geq \frac{2}{7}>0$.
If $n=1$, then $v$ is incident to at least one $5^+$-face, and therefore,
$ch^*(v)\geq \frac{2}{7}+\frac{1}{5}-\frac{1}{3}>0$ by rule R2. If $n=2$, then $v$ is incident to at least two $5^+$-faces,
and therefore $ch^*(v)\geq \frac{2}{7}+2\times\frac{1}{5}-2\times\frac{1}{3}=\frac{2}{105}>0$, by rule R2.
Hence, $0\leq\sum_{v\in V(G)}\frac{2}{7} + \sum_{f\in F(G)}(d_G(f)-4) \leq -8$, a contradiction.


\begin{thebibliography}{99}

\bibitem{Gruenewald_2000} S.~Gr\"unewald,~Chromatic Index Critical Graphs and Multigraphs, Dissertation, Fakult\"at f\"ur Mathematik, Universit\"at Bielefeld (2000).

\bibitem{Luo_Miao_Zhao_2009}  R. Luo, L. Miao and Y. Zhao,  The size of edge chromatic critical graphs with maximum degree 6, J.~Graph Theory {\bf 60} (2009) 149 - 171.

\bibitem{Sanders_Zhao_2001} D.~P.~Sanders, Y.~Zhao, Planar graphs of maximum degree seven are class 1, J.~Combin.~Theory Ser.~B {\bf 83} (2001) 201-212.

\bibitem{Vizing_1964} V.~G.~Vizing, On an estimate of the chromatic index of a p-graph, Metody Diskret.~Analiz {\bf 3} (1964) 25-30 (in Russian).

\bibitem{Wang_Xu_2013}  Y.~Wang, L.~Xu,  A sufficient condition for a plane graph with maximum degree 6 to be class 1, Discrete Appl. Math. {\bf 161} (2013) 307-310.

\bibitem{Woodall_2008} D.~R.~Woodall, The average degree of an edge-chromatic critical graph, Discrete Math.~{\bf 308} (2008) 803-819.

\bibitem{Zhang_2000} L.~Zhang, Every graph with maximum degree 7 is of class 1, Graphs Combin.~{\bf 16} (2000) 467-495.

\end{thebibliography}
\end{document}